\documentclass[10pt]{article}
\usepackage{amsmath,amsthm,amsfonts,amssymb,latexsym,graphicx,etoolbox,stmaryrd,environ}

\usepackage[shortcuts]{extdash}

\allowdisplaybreaks[4]

\usepackage[pdfpagemode=UseNone,pdfstartview=FitH]{hyperref}

\newtheorem{theorem}{Theorem}
\newtheorem{proposition}[theorem]{Proposition}

\theoremstyle{definition}

\newtheorem{example}[theorem]{Example}

\theoremstyle{remark}

\newcommand{\dd}{\,\mathrm{d}}

\renewcommand{\P}{\mathbb{P}}
\newcommand{\E}{\mathbb{E}}
\newcommand{\FFF}{\mathcal{F}}

\newcommand{\KL}{\textrm{KL}}

\DeclareMathOperator{\UEP}{\mathbb{P}^{\mathrm{exch}}}
\DeclareMathOperator{\UCP}{\mathbb{P}^{\mathrm{ce}}}

\title{Conformal e-testing}
\author{Vladimir Vovk, Ilia Nouretdinov, and Alex Gammerman}

\begin{document}
\maketitle

\begin{abstract}
  There is a useful counterpart of conformal prediction for e-values,
  called \emph{conformal e-prediction}.
  Conformal prediction can serve as basis for testing the assumption of exchangeability,
  leading to \emph{conformal testing}.
  Similarly, conformal e-prediction can also serve as basis for testing.
  The resulting \emph{conformal e-testing}
  looks very different from but inherits some strengths of conformal testing;
  it even has some advantages over conformal testing.
  In this paper we discuss systematically both strengths and limitations of conformal e-testing.

   The version of this paper at \url{http://alrw.net} (Working Paper 29)
   is updated most often.
\end{abstract}

\section{Introduction}

A useful application of conformal prediction is conformal testing,
which is a technique for testing the assumption of exchangeability
(or another online compression model).
Conformal e-prediction is a modification of conformal prediction
obtained by replacing the notion of a p-value by that of an e-value;
it is reviewed in the sister article \cite{Vovk:sister}.
(It may be natural to refer to conformal testing as \emph{conformal p-testing},
but we will never use this term.)

An important advantage of conformal prediction over conformal e-prediction
is that its strong property of validity found in \cite[Theorem~1]{Vovk:2002FOCS}
allows us to test the assumption of exchangeability.
This strong property can be stated as the independence of the smoothed conformal p-values
output at different steps,
and to test exchangeability,
we can bet against the conformal p-values being independent and uniformly distributed.
This led to the introduction in 2003 \cite{Vovk/etal:2003ICML}
of conformal test martingales.
The validity of conformal test martingales as means of testing exchangeability
shows, e.g., in Ville's theorem \cite[p.~100]{Ville:1939}:
under exchangeability, the probability that a given conformal test martingale $S$
ever exceeds a fixed threshold $c>1$ is at most $1/c$.
For example, we might feel justified in rejecting the hypothesis of exchangeability
when $S$ exceeds $100$,
since the probability of this event is at most $1\%$.

In general and informally, validity is the requirement that our testing methods
should give false evidence against exchangeability (in the context of this paper)
only with low probability.
It will appear in various guises in this paper.
Under the restriction of validity,
we would also like our procedures to be efficient
at discovering evidence against exchangeability.
For a long time, nothing was known about the efficiency of conformal test martingales,
and first results about their efficiency appeared in 2019
(see \cite{Vovk:2021-local});
this is a major topic of \cite[Part~III]{Vovk/etal:2022book}.

In this paper we discuss conformal e-testing systematically
and compare it with conformal testing.
It turns out that, similarly to the case of conformal e-prediction \cite{Vovk:sister},
conformal e-testing can often emulate strengths of conformal testing.
Moreover, conformal e-testing has some advantages of its own over conformal testing.

We start the main part of the paper in Sect.~\ref{sec:CeP}
by defining conformal e-prediction in a way adapted to the use in conformal e-testing.
In the following Sect.~\ref{sec:online} we introduce conformal e-testing and explore its validity.
Instead of conformal test martingales,
we obtain what we call ``conformal e-pseudomartingales'',
and our main finding here is negative:
conformal e-pseudomartingales can violate badly
the property of validity expressed by Ville's inequality.
While in Sect.~\ref{sec:online} we use the online setting, which is standard in this area,
in Sect.~\ref{sec:batch} we use a more limited batch setting,
which allows us to establish results about both validity and efficiency of conformal e-testing.

In Sections~\ref{sec:CUSUM-validity} and~\ref{sec:CUSUM-efficiency}
we discuss ``multistage'' ways of testing the exchangeability assumption,
concentrating on the standard CUSUM procedure used on top of conformal e-prediction,
which we call the ``conformal CUSUM e-procedure''.
Section~\ref{sec:CUSUM-validity} is devoted to the validity of the conformal CUSUM e-procedure:
we show that under exchangeability it raises false alarms with a frequency
determined by its parameter (and it can be made as low as we want).
This property is deduced from the analogous property of what we call the reverse Shiryaev--Roberts procedure.
In this section we also give an example showing an undesirable property of the Shiryaev--Roberts procedure,
standard or reverse.
In Sect.~\ref{sec:CUSUM-efficiency} we discuss a property of efficiency of the conformal CUSUM e-procedure.
While Sect.~\ref{sec:CUSUM-validity} shows that conformal e-testing inherits some strengths of conformal testing,
Sect.~\ref{sec:CUSUM-efficiency} demonstrates an advantage of conformal e-testing.

Section~\ref{sec:conclusion} concludes and lists some directions of further research.

\section{Conformal e-prediction}
\label{sec:CeP}

The task of conformal e-prediction, and predictive machine learning in general,
is to predict the label of a test object $x$
given a training set $z_1,\dots,z_n$ whose elements are labelled objects $z_i=(x_i,y_i)$.
What distinguishes conformal (e-) prediction is that
for each potential label $y$ for the test object $x$ it provides a nonnegative number
$f(z_1,\dots,z_n,(x,y))$
(we usually drop the internal parentheses)
reflecting the plausibility of $y$ being the true label of $x$.

The objects $x_i$ are drawn from the \emph{object space} $\mathbf{X}$
and the labels $y_i$ from the \emph{label space} $\mathbf{Y}$;
both are required to be non-empty measurable spaces.
The \emph{observations} $z=(x,y)$ are drawn from the Cartesian product
(the \emph{observation space}) $\mathbf{Z}:=\mathbf{X}\times\mathbf{Y}$.
In this paper we will also be interested in the case where $\mathbf{Z}$ is unstructured
(not a Cartesian product, which can be embedded into the structured case
by setting $\mathbf{X}$ or $\mathbf{Y}$ to a one-element space).

We will use the notation $X^+:=\cup_{n=1}^{\infty}X^n$
for the set of all non-empty finite sequences of elements of $X$.
If $X$ is a measurable space, $X^+$ is also a measurable space.

A \emph{nonconformity e-measure} is a measurable function $A:\mathbf{Z}^+\to[0,\infty)^+$
that maps every finite sequence $(z_1,\dots,z_m)$, $m\in\{1,2,\dots\}$,
to a finite sequence $(\alpha_1,\dots,\alpha_m)$ of the same length such that
\begin{equation}\label{eq:1}
  \frac1m
  \sum_{i=1}^m
  \alpha_i
  \le
  1
\end{equation}
and that satisfies the property of \emph{equivariance}:
for any $m$ and any permutation $\pi$ of $\{1,\dots,m\}$,
\[
  (\alpha_1,\dots,\alpha_m) = A(z_1,\dots,z_m)
  \Longrightarrow
  (\alpha_{\pi(1)},\dots,\alpha_{\pi(m)}) = A(z_{\pi(1)},\dots,z_{\pi(m)}).
\]
(We sometimes refer to the $\alpha_i$ as \emph{nonconformity e-scores}.)
The corresponding \emph{conformal e-predictor}
$f:\mathbf{Z}^+\to[0,\infty)$
is defined as
\begin{equation*}
  f(z_1,\dots,z_n,x,y)
  :=
  \alpha_{n+1},
  \text{\enspace where\enspace}
  (\alpha_1,\dots,\alpha_n,\alpha_{n+1}) := A(z_1,\dots,z_n,(x,y)).
\end{equation*}

For a training set $z_1,\dots,z_n$ and a test object $x$,
the full prediction for $x$ according to a conformal e-predictor $f$
is given by the family of \emph{potential conformal e-values}
\[
  \left(
    f(z_1,\dots,z_n,x,y)
    \mid
    y\in\mathbf{Y}
  \right).
\]
We can make a confident prediction for $x$ if the potential conformal e-values are large
for all $y\in\mathbf{Y}$ except for one.

A nonconformity e-measure and the corresponding conformal e-predictor
are \emph{admissible}
if we always have ``$=$'' in place of ``$\le$'' in~\eqref{eq:1} in its definition.
Testing procedures based on nonconformity e-measures that are not admissible
can be improved.

Now let us state what we regard as the main property of validity of conformal e-prediction
(for more information, see \cite[Sect.~3]{Vovk:sister}).
Let $Z_1,Z_2,\dots$ be the random observations,
i.e., the random elements whose realizations are the observed $z_1,z_2,\dots$.
In general, $(X,Y)$ or $Z$ are random observations,
i.e., random elements taking values in the observation space $\mathbf{Z}$.
Remember that a finite sequence of random elements is \emph{exchangeable}
if its joint distribution does not change if it is permuted
(and an infinite sequence is exchangeable
if its every finite beginning is exchangeable).
For example, any IID sequence
(i.e., a sequence of independent and identically distributed random elements)
is exchangeable.

With each conformal e-predictor $f$ we can associate the sequence of \emph{conformal e-variables}
\begin{equation}\label{eq:E}
  E_n:=f(Z_1,\dots,Z_{n-1},Z_n).
\end{equation}
Intuitively, large values of the conformal e-variables are evidence
against $Z_1,Z_2,\dots$ being exchangeable (in particular, IID).

The \emph{exchangeable filtration} \cite[Sect.~5.6]{Ramdas/etal:2023} is $(\FFF_n)$,
where $\FFF_n$ is the $\sigma$-algebra generated by the multiset $\lbag Z_1,\dots,Z_{n-1}\rbag$
and the observations $Z_n,Z_{n+1},\dots$.
The following proposition is a stronger version of the property of validity of conformal e-prediction.

\begin{proposition}\label{prop:validity-1}
  For any conformal e-predictor $f$ and any $n$,
  if the sequence $Z_1,Z_2,\dots$ is exchangeable,
  then
  \begin{equation}\label{eq:validity}
    \E(E_n\mid\FFF_{n+1})\le1,
  \end{equation}
  where $\FFF$ is the exchangeable filtration and
  $E_n$ is the conformal e\-/variable~\eqref{eq:E}
  (with ``$=$'' in place of ``$\le$'' in \eqref{eq:validity}
  if $f$ is admissible).
\end{proposition}

\begin{proof}
  By the definition of conformal e-predictors we have
  \[
    \E
    \left(
      f(Z_1,\dots,Z_n)
      \mid
      \lbag Z_1,\dots,Z_n\rbag
    \right)
    \le
    1,
  \]
  and we can add $Z_{n+1},Z_{n+2},\dots$ to the condition
  since $Z_1,\dots,Z_n$ are exchangeable conditionally on $Z_{n+1},Z_{n+2},\dots$. 
  This is equivalent to \eqref{eq:validity}.
\end{proof}

We still have \eqref{eq:validity} if the sequence $Z_1,Z_2,\ldots = Z_1,\dots,Z_N$ is finite
provided its length $N$ is at least $n$, $N\ge n$.
In this case the exchangeable filtration is $(\FFF_n)_{n\le N+1}$,
where $\FFF_{N+1}$ is the $\sigma$-algebra generated by the multiset $\lbag Z_1,\dots,Z_N\rbag$.

A weaker version of \eqref{eq:validity} is $\E(E_n)\le1$.
This can be expressed as $E_n$ being an e-variable,
where an \emph{e-variable} is defined to be a nonnegative random variable
with expected value at most 1.
The values taken by e-variables are referred to as \emph{e-values}.

\section{Conformal e-testing in the online setting}
\label{sec:online}

In this section we start our discussion of testing the assumption of exchangeability
from the online setting.
Namely, the testing process proceeds in time by processing
a potentially infinite stream of observations $z_1,z_2,\dots$
sequentially one by one,
and at each moment we would like to have a measure of the amount of evidence
that we have found against the assumption of their exchangeability.
If such a measure exceeds a large threshold $c$, we might want to raise an alarm
indicating that exchangeability is likely to have been violated.
In conformal testing, such a measure is provided by conformal test martingales,
and it satisfies a strong requirement of validity
(an application of Ville's inequality):
we will raise an alarm with probability at most $1/c$.
In this section we will see that the natural counterpart of conformal test martingales
in conformal e-testing,
which we call conformal e-pseudomartingales,
fails completely to satisfy this strong property of validity
(Proposition~\ref{prop:fail} below).

The \emph{conformal e-pseudomartingale}
corresponding to the conformal e-variables \eqref{eq:E} is
\begin{equation*}
  S_n
  :=
  E_1 \dots E_n,
  \quad
  n=0,1,2,\dots,
\end{equation*}
where $S_0$ is understood to be 1.
It may not be a genuine martingale
since by Proposition~\ref{prop:validity-1} we have $\E(E_n\mid\FFF_{n+1})=1$ for all $n$
instead of $\E(E_n\mid E_1,\dots,E_{n-1})=1$ required in the definition of martingales.

The definition of conformal e-pseudomartingales, nevertheless,
is very similar to that of conformal test martingales,
and the conformal e-variables $E_n$ look analogous to the betting functions
of conformal testing \cite[Sect.~8.1.2]{Vovk/etal:2022book}.
It corresponds to the gambling picture
in which we start from an initial capital of $1$
and then compound the conformal e-values as usual by multiplying them.
The crucial difference is that in conformal testing
a betting function only depends on the past p-values,
whereas in conformal e-testing it may also depend on the multiset of actual observations.

The following proposition uses the notation
\[
  S^*_{\infty}
  :=
  \sup_{n=1,2,\dots}
  S_n.
\]
Ville's inequality can then be written as $P(S^*_{\infty}\ge c)\le1/c$
for any $c>1$ and any test martingale $S$ w.r.\ to $P$
(i.e., any nonnegative martingale $S$ satisfying $S_0=1$).
To exclude the trivial case,
let us assume that the $\sigma$-algebra on $\mathbf{Z}$ is different from $\{\emptyset,\mathbf{Z}\}$;
in particular, $\left|\mathbf{Z}\right|>1$.

\begin{proposition}\label{prop:fail}
  For any $\epsilon>0$ and $c>1$,
  there exists an exchangeable probability measure $P$ on $\mathbf{Z}^{\infty}$
  and a conformal e-pseudomartingale such that $P(S^*_{\infty}\ge c)\ge1-\epsilon$.
\end{proposition}

\begin{proof}
  Assume, without loss of generality, that $\mathbf{Z}=\{0,1\}$
  and that $c$ is an integer.
  Consider an e-predictor $f$ that stakes everything on 1;
  in particular, for sequences of any length $n$,
  \[
    f:(0,\dots,0,1)\mapsto (0,\dots,0,n);
  \]
  on the other hand, let
  \[
    f:(0,\dots,0)\mapsto (1,\dots,1).
  \]
  Let $A_c\subseteq\{0,1\}^{\infty}$ be the set of all sequences
  that have $0^{[c]}$ as their prefix apart from the sequence $0^{[\infty]}$
  (where $0^{[a]}$ is the sequence consisting of $a$ 0s).
  Then $A_c\subseteq\{S^*_{\infty}\ge c\}$;
  besides, $P(A_c)\ge1-\epsilon$ if, under $P$,
  the random observations $Z_1,Z_2,\dots$ are generated in the IID manner
  with the probability of 1 sufficiently small (but positive).
\end{proof}

We interpret the statement of Proposition~\ref{prop:fail}
as a complete loss of validity (in the strong sense of Ville's inequality)
for conformal e-pseudomartingales.
In particular, Proposition~\ref{prop:fail}
demonstrates that conformal e-pseudomartingales are not martingales in general.
Despite conformal e-pseudomartingales violating Ville's inequality so badly,
we will see in Sect.~\ref{sec:CUSUM-validity}
that the CUSUM procedure based on conformal e-pseudomartingales
still satisfies the standard property of validity
(as in \cite[Corollary~8.14]{Vovk/etal:2022book}).
Besides, the strong property of validity holds for a fixed time horizon,
as Proposition~\ref{prop:validity-2} below will show.

\section{Conformal e-testing in the batch setting}
\label{sec:batch}

Having established loss of validity in the online setting,
in this section we move on to the batch setting,
which is standard in statistics
(although we use terminology that is standard in conformal testing rather than statistics).
Namely, we assume that the number $N$ of observations is fixed and known in advance.
At the end of step $N$ we are required to make a decision
whether to reject the hypothesis of exchangeability,
or at least to present the amount of evidence that we have found
against the null hypothesis of exchangeability.

Conformal e-testing in the batch mode is defined in two steps.
First, a function $E:\mathbf{Z}^N\to[0,\infty)$
is called a \emph{basic conformal e-test} if it has the form
\[
  E(z_1,\dots,z_N)
  =
  \prod_{n=1}^N
  f(z_1,\dots,z_n),
  \quad
  \forall(z_1,\dots,z_N)\in\mathbf{Z}^N,
\]
for some conformal e-predictor $f$.
In other words, if $E=S_N$ for some conformal e-pseudomartingale $S$.
And second, a \emph{conformal e-test} is defined
as a convex combination of basic conformal e-tests:
namely, the conformal e-tests are defined as convex mixtures
$\lambda_1 E^{(1)}+\dots+\lambda_k E^{(k)}$ of basic conformal e-tests $E^{(1)},\dots,E^{(k)}$,
where $k\in\{1,2,\dots\}$, $\lambda_i\in[0,1]$, and $\lambda_1+\dots+\lambda_k=1$.
We need the second step to achieve some efficiency later on,
and the following proposition gives a property of validity
for conformal e-pseudomartingales at a fixed time
and, slightly more generally, for conformal e-tests.

\begin{proposition}\label{prop:validity-2}
  Suppose $Z_1,\dots,Z_N$ are exchangeable.
  For any conformal e-test $E$,
  $E(Z_1,\dots,Z_N)$ is a bona fide e-variable.
  It satisfies $\E(E(Z_1,\dots,Z_N))=1$ if the underlying conformal e-predictor is admissible.
\end{proposition}

\begin{proof}
  Assume, without loss of generality, that $E$ is a basic conformal e-test.
  We will use the notation \eqref{eq:E},
  where $f$ is the underlying conformal e-predictor,
  and the exchangeable filtration $(\FFF_n)$.
  We can show that $\E(E_1\dots E_n\mid\FFF_{n+1})\le1$ a.s.\ by induction in $n$
  (for $n=0$ this statement is vacuously true as equality):
  \begin{multline*}
    \E(E_1\dots E_n\mid\FFF_{n+1})
    =
    \E(\E(E_1\dots E_n\mid\FFF_{n})\mid\FFF_{n+1})\\
    =
    \E(E_n\E(E_1\dots E_{n-1}\mid\FFF_{n})\mid\FFF_{n+1})
    \le
    \E(E_n\mid\FFF_{n+1})
    \le
    1
    \quad\text{a.s.},
  \end{multline*}
  the last inequality following from Proposition~\ref{prop:validity-1}.
  If the underlying conformal e-predictor is admissible,
  both ``$\le$'' become ``$=$''.
\end{proof}

If a conformal e-test $E$ is chosen in advance and takes a very large value on the realized data sequence,
we are justified in rejecting the assumption of exchangeability.

As the next step,
we explore the efficiency of conformal e-testing
along the lines of the treatment of the efficiency of conformal testing
in \cite[Sect.~9.1]{Vovk/etal:2022book}.
As in the case of conformal testing \cite[Sect.~9.1]{Vovk/etal:2022book},
for our (rather weak) statement of efficiency
we will simplify our task by only considering the binary case,
$\mathbf{Z}=\{0,1\}$.

First we define unrestricted testing of exchangeability.
The \emph{upper exchangeability probability} $\UEP$ of measurable $A\subseteq\mathbf{Z}^N$
is defined as
\begin{equation}\label{eq:UEP}
  \UEP(A)
  :=
  \sup_P
  P(A),
\end{equation}
$P$ ranging over the exchangeable probability measures on $\mathbf{Z}^N$.
The intuition behind $\UEP(A)$ is that,
if it is very small,
$A$ can be used for testing the exchangeability of the data-generating distribution:
we are entitled to reject exchangeability if $(z_1,\dots,z_N)\in A$
provided $A$ is chosen in advance.

The \emph{upper conformal e-probability} $\UCP$ of $A\subseteq\mathbf{Z}^N$
is defined as
\begin{equation}\label{eq:UCP}
  \UCP(A)
  :=
  \inf
  \left\{
    \alpha:
    \exists E
    \;
    \forall(z_1,\dots,z_N)\in A:
    E(z_1,\dots,z_N)
    \ge
    1/\alpha
  \right\},
\end{equation}
$E$ ranging over the conformal e-tests.
The intuition is the same:
if $\UCP(A)$ is very small,
we are entitled to reject exchangeability if we observe $(z_1,\dots,z_N)\in A$,
again assuming that $A$ is chosen in advance.
The difference from $\UEP(A)$ is that now the lack of exchangeability
must be demonstrated via conformal e-testing.

The following proposition is an analogue of \cite[Proposition~9.5]{Vovk/etal:2022book}.
Intuitively, its second statement says that conformal e-testing is universal
in the batch mode:
if lack of exchangeability can be demonstrated at all,
it can be demonstrated (albeit less convincingly) using conformal e-testing.
(And its first statement is another expression of validity.)

\begin{proposition}
  For any event $A\subseteq\mathbf{Z}^N$,
  $\UEP(A)\le\UCP(A)$.
  Assuming $\mathbf{Z}=\{0,1\}$,
  $\UCP(A)\le N\UEP(A)$.
\end{proposition}

\begin{proof}
  The inequality $\UEP(A)\le\UCP(A)$ follows immediately from Markov's inequality
  applied to conformal e-tests $E$ in combination with Proposition~\ref{prop:validity-2}:
  if $E(z_1,\dots,z_N)\ge1/\alpha$ for all $(z_1,\dots,z_N)\in A$,
  \[
    \UEP(A)
    \le
    \P(E(Z_1,\dots,Z_N)\ge1/\alpha)
    \le
    \frac{\E(E(Z_1,\dots,Z_N))}{1/\alpha}
    \le
    \alpha,
  \]
  assuming $Z_1,\dots,Z_N$ are exchangeable.

  Now assume $\mathbf{Z}=\{0,1\}$;
  this part of the proof will be a modification of the proof of Proposition~9.5
  in \cite[Sect.~9.4.2]{Vovk/etal:2022book}.
  As a first step, notice that it suffices to prove $\UCP(A)\le\UEP(A)$
  for any nonempty $A\subseteq\mathbf{Z}^N$
  such that each sequence in $A$ has the same number of 1s.
  Let us fix such an $A$, and let $K\in\{0,\dots,N\}$ be the number of 1s
  in the elements of $A$.
  For each sequence $\zeta=(z_1,\dots,z_N)\in A$,
  consider the basic conformal e-test $E_{\zeta}=e_1\dots e_N$,
  where the $n$th e-value $e_n$ is
  \[
    e_n
    :=
    \begin{cases}
      n/k & \text{if $z_n=1$}\\
      n/(n-k) & \text{if $z_n=0$},
    \end{cases}
  \]
  $k$ being the number of 1s among the first $n$ elements of $\zeta$.
  (This corresponds to a nonconformity e-measure satisfying
  \[
    A
    \left(
      0^{[n-k]},1^{[k]}
    \right)
    :=
    \begin{cases}
      \left(0^{[n-k]},\left(\frac{n}{k}\right)^{[k]}\right) & \text{if $z_n=1$}\\
      \left(\left(\frac{n}{n-k}\right)^{[n-k]},0^{[k]}\right) & \text{if $z_n=0$},
    \end{cases}
  \]
  where $b^{[a]}=b,\dots,b$ ($a$ times), as in the proof of Proposition~\ref{prop:fail}.
  This was called ``reckless gambling'' in \cite[end of Sect.~9.1.2]{Vovk/etal:2022book}.)
  The product $E_{\zeta}=e_1\dots e_N$ will then be
  \[
    \frac{N!}{K!(N-K)!}=\binom{N}{K}.
  \]
  The arithmetic mean of $E_{\zeta}$ over $\zeta\in A$ witnesses that
  \[
    \UCP(A)
    \le
    \lvert A\rvert / \binom{N}{K}
    =
    \UEP(A).
    \qedhere
  \]
\end{proof}

Proposition~\ref{prop:fail} can be restated in terms of $\UEP$ and $\UCP$
adapted to the online setting.
Let $A\subseteq\mathbf{Z}^{\infty}$.
Define $\UEP(A)$ by~\eqref{eq:UEP}, as before,
with $P$ ranging over the exchangeable probability measures on $\mathbf{Z}^{\infty}$.
In the spirit of Ville's inequality,
let us modify \eqref{eq:UCP} as
\begin{equation*}
  \UCP(A)
  :=
  \inf
  \left\{
    \alpha:
    \exists S
    \;
    \forall(z_1,z_2,\dots)\in A:
    S^*_{\infty}(z_1,z_2,\dots)
    \ge
    1/\alpha
  \right\},
\end{equation*}
where $S$ ranges over the conformal e-pseudomartingales
and we slightly abuse our notation by regarding $S_n$ as functions of the observations.
Then we can see that $\UEP(A_c)=1$ and $\UCP(A_c)\le1/c\to0$ as $c\to\infty$
(where $A_c$ are defined in the proof of Proposition~\ref{prop:fail}).
This can be expressed as $\UEP$ and $\UCP$ being entirely asymptotically separated
(cf.~\cite[Sect.~3.10, especially (3)]{Shiryaev:2016}).

\section{The conformal CUSUM e-procedure: validity}
\label{sec:CUSUM-validity}

In the previous section we discussed testing the assumption of exchangeability once,
but in some important applications we would like to test it repeatedly over time
in the online mode (see, e.g., \cite[Sect.~8.3]{Vovk/etal:2022book}):
as soon as we suspect that exchangeability is violated, we raise an alarm,
and we are allowed do so more than once.
The topic of this and next sections is such multistage exchangeability testing.
As usual, such procedures are required to satisfy
properties of validity and efficiency.
Here validity means that, under exchangeability,
the probability or frequency of alarms
(which in this case are \emph{false alarms})
should be bounded by a prespecified constant.
When applied to the moment when the first alarm is raised,
such properties of validity are much weaker
than the property of validity discussed in the previous sections:
the first alarm is usually raised, sooner or later, with probability one.
And efficiency means that, if exchangeability is violated at some point,
an alarm should be raised as quickly as possible afterwards.
In this section we concentrate on validity of multistage testing.

Suppose we observe a sequence of e-values $e_1,e_2,\dots$
output by a conformal e-predictor in the online protocol,
as described in Sect.~\ref{sec:online},
under exchangeability.
The \emph{conformal CUSUM e-procedure} \cite{Page:1954} raises the $k$th alarm,
$k=1,2,\dots$,
at the time
\begin{equation}\label{eq:e-CUSUM}
  \tau_k
  :=
  \min
  \left\{
    n>\tau_{k-1}:
    \max_{i\in\{\tau_{k-1}+1,\dots,n\}}
    e_i\dots e_n
    \ge
    c
  \right\},
\end{equation}
where $\tau_0:=0$ and $c>1$ is the parameter of the procedure.
It is usually applied in the situation where the observations $z_i$
are generated independently
first from a known probability measure $Q_0\in\mathfrak{P}(\mathbf{Z})$
and then from another known probability measure $Q_1\in\mathfrak{P}(\mathbf{Z})$
($\mathfrak{P}(\mathbf{Z})$ being the family of all probability measures on $\mathbf{Z}$),
and where $e_i$ is the likelihood ratio of $Q_1$ to $Q_0$ evaluated at $z_i$
(see, e.g., \cite[Sect.~6.2]{Poor/Hadjiliadis:2009} or \eqref{eq:L} below).
In our current context, where the $e_i$ are conformal e-values,
we may call it the \emph{conformal CUSUM e\-/procedure}.

The following proposition
gives an asymptotic property of validity
for the conformal CUSUM e-procedure.

\begin{proposition}\label{prop:CUSUM}
  Let $A_n$ be the number of alarms
  \begin{equation}\label{eq:A}
    A_n
    :=
    \max\{k\mid\tau_k\le n\}
  \end{equation}
  raised by the conformal CUSUM e-procedure \eqref{eq:e-CUSUM}
  after processing the first $n$ observations $Z_1,\dots,Z_n$.
  Then
  \begin{equation}\label{eq:SLLN}
    \limsup_{n\to\infty}
    \frac{A_n}{n}
    \le
    \frac1c
    \quad
    \text{a.s.}
  \end{equation}
  provided the observations $Z_1,Z_2,\dots$ are exchangeable.
\end{proposition}

In practical applications,
Proposition~\ref{prop:CUSUM} can be applied when,
after investigating each alarm, it is decided that the alarm was false,
and so we can continue processing a single stream of observations.
If the alarm was genuine,
we need to reset the multistage procedure,
and Proposition~\ref{prop:CUSUM} is not applicable.

A popular modification of the CUSUM procedure
is the \emph{Shiryaev--Roberts procedure},
which replaces the $\max$ in \eqref{eq:e-CUSUM}
by the sum $\sum$ (see, e.g., \cite[Sect.~6.5]{Poor/Hadjiliadis:2009}).
To state our next validity result (Proposition~\ref{prop:rSR} below),
we need the modification
\begin{equation}\label{eq:rSR}
  \tau_k
  :=
  \min
  \left\{
    n>\tau_{k-1}:
    \max_{i\in\{\tau_{k-1}+1,\dots,n\}}
    \sum_{j=i}^n
    e_i\dots e_j
    \ge
    c
  \right\},
\end{equation}
which we will call the \emph{reverse Shiryaev--Roberts procedure}
(or \emph{conformal reverse Shiryaev--Roberts e-procedure}
to emphasize $e_i$ being produced by a conformal e-predictor).

\begin{proposition}\label{prop:rSR}
  Let $A_n$ be defined by~\eqref{eq:A}
  for the conformal reverse Shiryaev--Roberts e-procedure \eqref{eq:rSR}.
  Then we still have~\eqref{eq:SLLN}
  provided the observations $Z_1,Z_2,\dots$ are exchangeable.
\end{proposition}

See Appendix~\ref{app:proofs-validity} for the proofs.
We will, of course, deduce Proposition~\ref{prop:CUSUM}
from Proposition~\ref{prop:rSR}:
it is clear that the Shiryaev--Roberts procedure
(regular or reverse)
raises alarms at least as often as CUSUM does.
However, the following example shows an advantage of CUSUM
from an intuitive point of view.

\begin{example}\label{ex:SR}
  Consider the vacuous admissible conformal e-predictor identically equal to 1.
  The CUSUM procedure based on it will never raise alarms,
  while the Shiryaev--Roberts e-procedure will raise alarms every $\lceil c\rceil$th step,
  thereby fully exploiting (for an integer $c$) the leeway
  permitted by our target property of validity~\eqref{eq:SLLN}.
  This example shows one feature of the Shiryaev--Roberts procedure
  (shared by the reverse Shiryaev--Roberts procedure)
  that can be considered its disadvantage:
  while the CUSUM procedure raises an alarm when it has genuine evidence for disorder,
  the Shiryaev--Roberts procedure may raise an alarm
  simply because it is allowed to do so by a given constraint;
  we might not have any evidence for disorder.
\end{example}

Because of the feature of the Shiryaev--Roberts procedure
illustrated in Example~\ref{ex:SR},
in the next section we will concentrate on CUSUM-type procedures.

\section{The conformal CUSUM e-procedure: efficiency}
\label{sec:CUSUM-efficiency}

The standard version of the CUSUM procedure satisfies
important properties of optimality \cite[Chap.~6]{Poor/Hadjiliadis:2009}:
it is optimal in Lorden's \cite{Lorden:1971} sense,
as shown by Moustakides \cite{Moustakides:1986};
it is also optimal in Ritov's \cite{Ritov:1990} very natural game-theoretic sense.
However, as we mentioned in the previous section,
the standard CUSUM procedure works under restrictive assumptions:
we know the prechange distribution $Q_0\in\mathfrak{P}(\mathbf{Z})$
and the postchange distribution $Q_1\in\mathfrak{P}(\mathbf{Z})$,
and the only unknown is the changepoint $N_0$;
the observations $Z_1,\dots,Z_{N_0}$ are generated from $Q_0$,
and the observations $Z_{N_0+1},Z_{N_0+2},\dots$ are generated from $Q_1$,
all independently.
Let $f_0$ and $f_1$ be probability densities of $Q_0$ and $Q_1$, respectively,
w.r.\ to a $\sigma$-finite measure $\mu$ on $\mathbf{Z}$ (such as $\mu:=Q_0+Q_1$)
under which both $Q_0$ and $Q_1$ are absolutely continuous.
We will use $\E_0$ and $\E_1$ for expectations w.r.\ to $Q_0$ and $Q_1$, respectively.

The \emph{standard CUSUM procedure} is based on the likelihood ratios
\begin{equation}\label{eq:L}
  L_n
  :=
  f_1(z_n) / f_0(z_n).
\end{equation}
(For simplicity, the reader may assume that $Q_0$ and $Q_1$
are positive discrete distributions
and replace $f_0(z)$ and $f_1(z)$ by $Q_0(\{z\})$ and $Q_1(\{z\})$,
respectively.)
As in \eqref{eq:e-CUSUM},
the $k$th alarm is raised at the time
\begin{equation}\label{eq:CUSUM}
  \tau_k
  :=
  \min
  \left\{
    n>\tau_{k-1}:
    \max_{i\in\{\tau_{k-1}+1,\dots,n\}}
    L_i\dots L_n
    \ge
    c
  \right\},
\end{equation}
with $\tau_0:=0$.
For comparison, we will also discuss experimental results
for the \emph{conformal CUSUM procedure},
which is defined in a similar way using conformal p-values,
as explained in \cite[Sect.~8.3.1]{Vovk/etal:2022book}.

This section implements a version of the ``Burnaev--Wasserman programme''
\cite[Sect.~2.5]{Vovk/etal:2022book} applied to the CUSUM procedure.
Suppose we would like to detect deviations from exchangeability,
but we have a prior model of the data generation mechanism
in which the observations are first generated from a given $Q_0\in\mathfrak{P}(\mathbf{Z})$
and then from another given probability measure, $Q_1\in\mathfrak{P}(\mathbf{Z})$,
independently of the previous observations.
We however, do not trust our model
and do not want the validity of our procedure to depend on it.
Therefore, we ``conformalize'' the standard CUSUM procedure,
using the likelihood ratios \eqref{eq:L}
(normalized to ensure \eqref{eq:1})
as nonconformity e-scores.
(Details will follow shortly.)
The procedure of conformalization may be said to work well
if the quality of the conformalized CUSUM is not significantly worse
than the quality of the standard CUSUM
even when the assumptions on which the standard CUSUM is based
are fully satisfied.

\begin{figure}[bt]
  \begin{center}
    \includegraphics[width=0.48\textwidth]{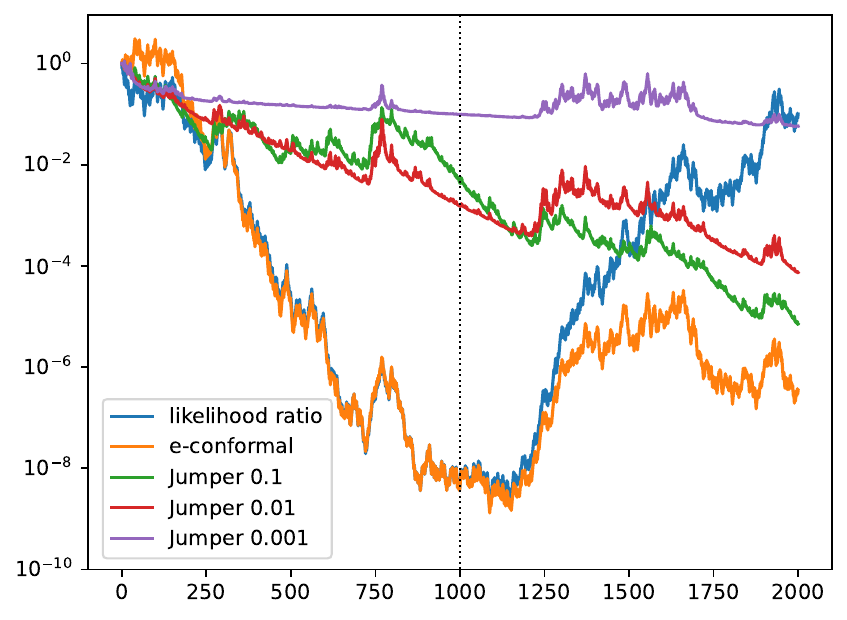}
    \includegraphics[width=0.48\textwidth]{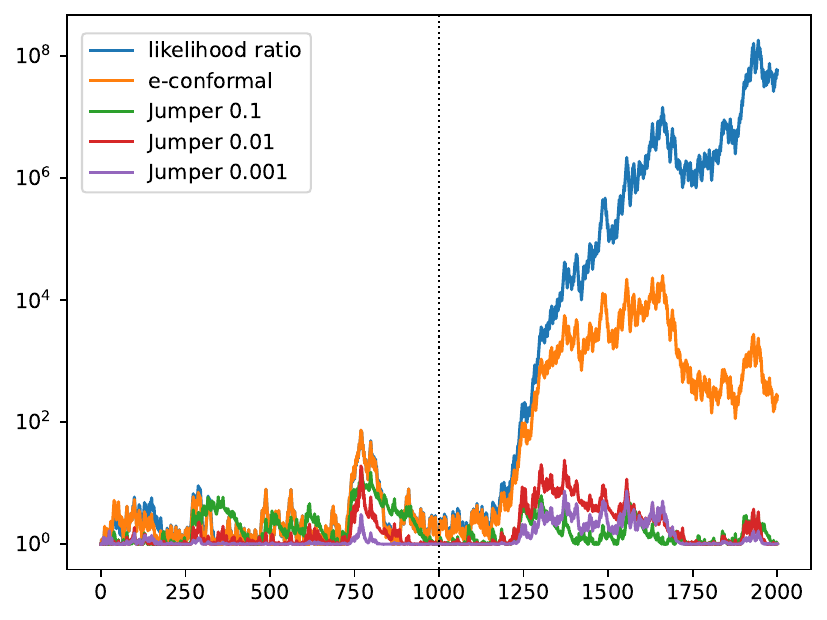}
  \end{center}
  \caption{Five stochastic processes, as described in text, in the Bernoulli case
    (with the parameter $0.5$ before the changepoint and $0.6$ after the changepoint).
    Left panel: the paths of the processes.
    Right panel: the paths of the corresponding CUSUM statistics.}
  \label{fig:Bernoulli}
\end{figure}

Some experimental results for the case of Bernoulli observations
are shown in Figure~\ref{fig:Bernoulli}.
In this experiment we generate $N_0:=1000$ observations $Z_1,\dots,Z_{N_0}$
from the Bernoulli distribution with parameter $0.5$
and another $N_1:=1000$ observations $Z_{N_0+1},\dots,Z_{N_0+N_1}$
from the Bernoulli distribution with parameter $0.6$.
The changepoint 1000 is shown as a dashed vertical line.
The left panel of Figure~\ref{fig:Bernoulli} shows the paths of five stochastic processes:
\begin{itemize}
\item
  the likelihood ratio martingale $L_1\dots L_n$, $n=0,1,\dots,N_0+N_1$, in blue;
\item
  the conformal e-pseudomartingale $E_1\dots E_n$, $n=0,1,\dots,N_0+N_1$ in orange;
  the conformal e-values are defined as the normalized $L_n$:
  \begin{equation}\label{eq:E-from-L}
    E_n
    :=
    \frac{L_n}{\frac1n(L_1+\dots+L_n)};
  \end{equation}
\item
  the green, red, and purple lines correspond to conformal testing;
  they will be discussed later and can be ignored for now.
\end{itemize}
The right panel of Figure~\ref{fig:Bernoulli} shows the corresponding \emph{CUSUM statistics},
where the CUSUM statistics of a process $S_n$, $n=0,1,2,\dots$ are defined as
\[
  S'_n
  :=
  S_n/\min(S_0,\dots,S_{n-1})
\]
(notice that $\tau_1$ in \eqref{eq:e-CUSUM} and \eqref{eq:CUSUM} is defined
as the moment when the CUSUM statistic first reaches level $c$).

The five lines (CUSUM statistics) in the right panel of Figure~\ref{fig:Bernoulli}
illustrate the efficiency of various versions of the CUSUM procedure
in detecting the changepoint.
The detection happens when the line first reaches the level $c$
after the changepoint.
(We ignore the rare case where the level $c$ is reached
between the changepoint and the last time before the changepoint
when the CUSUM statistic is 0.)
Without fixing $c$, we can judge how efficient the procedure is
by the slope of the line after the changepoint.
For example, we can see that the conformal CUSUM e-procedure
(corresponding to the orange line)
will raise an alarm not much later than the standard CUSUM procedure
(corresponding to the blue line)
for $c$ up to about $10^3$.

The blue and orange lines in the right panel are close to each other
shortly after the changepoint,
suggesting that conformal e-testing is efficient at first,
but then they start diverging.
This divergence illustrates the phenomenon of ``decay'' discussed in \cite[Sect.~8.4.1]{Vovk/etal:2022book}.
In fact, decay is inevitable and by itself does not indicate lack of efficiency of conformal e-testing
(unless it sets in too early).
Remember that the likelihood ratio martingale and conformal e-pseudomartingale
are testing very different null hypotheses:
for the former the null hypothesis is $Q_0$, and for the latter it is exchangeability.
If our null hypothesis is $Q_0$, we will be constantly surprised seeing observations from $Q_1\ne Q_0$,
whereas under exchangeability $Q_1$ will gradually become ``the new normal'',
and we will stop being surprised.
Therefore, decay will inevitably happen,
and it is the moment when it sets in that is the hallmark of efficiency
(the later the better).

In our experiments we always use the standard seed 42 for the numpy random number generator,
but the results are qualitatively similar for other seeds as well
(one unusual feature of Figure~\ref{fig:Bernoulli} is that the likelihood ratio martingale
does not start its ascent right after the changepoint,
which typically happens for other seeds).

The green, red, and purple lines in the left panel of Figure~\ref{fig:Bernoulli}
show the paths of three conformal test martingales,
as defined in \cite[Sect.~8.1.2]{Vovk/etal:2022book}.
(This and next paragraphs depend on and use the terminology of \cite[Sect.~8.1]{Vovk/etal:2022book}.)
The nonconformity score of each observation $Z_n$
(here and in Figure~\ref{fig:Cauchy} below) is defined to be
the likelihood ratio $L_n$
(which is equivalent to using $Z_n$ itself as nonconformity score in the current Bernoulli case).
This produces a sequence of conformal p-values $p_1,p_2,\dots$.
To turn conformal p-values into a conformal test martingale,
we need to define betting functions.
In Figure~\ref{fig:Bernoulli} we use the Simple Jumper betting functions
\cite[Sect.~8.1.2]{Vovk/etal:2022book}
omitting the adjective ``Simple'' in the caption.
The resulting Simple Jumper martingale depends on a parameter called \emph{jumping rate},
and we use three values for it as indicated in the legend.

It would be more in the spirit of the Burnaev--Wasserman programme
to adapt the betting functions to the assumed $Q_0$ and $Q_1$.
This is an interesting direction for further research,
but in this paper we are only using Simple Jumper,
a generic betting martingale that often gives satisfactory results
for a wide range of datasets \cite[Sect.~8.1.2, (8.8)]{Vovk/etal:2022book}.
(Using a generic betting martingale is analogous to our use
of a generic way of betting, given by \eqref{eq:E-from-L},
in the conformal CUSUM e-procedure.)

To start theoretical analysis of the phenomena demonstrated by the blue and orange lines
in Figure~\ref{fig:Bernoulli},
let us see why the blue line changes its slope after the changepoint.
Before the changepoint, $n\le N_0$, the expectation of the likelihood ratio
is at most 1:
\begin{equation}\label{eq:before-changepoint}
  \E_0 L_n
  =
  \E_0\frac{f_1(Z_n)}{f_0(Z_n)}
  =
  \int \frac{f_1}{f_0}f_0 \dd\mu
  \le
  \int f_1 \dd\mu
  =
  1;
\end{equation}
it is 1 if $Q_1$ is absolutely continuous w.r.\ to $Q_0$
(which is always the case in our experiments).
After the changepoint, $n>N_0$, it becomes
\begin{equation}\label{eq:after-changepoint}
  \E_1 L_n
  =
  \E_1\frac{f_1(Z_n)}{f_0(Z_n)}
  =
  \int \frac{f_1}{f_0}f_1 \dd\mu
  =
  1 + \int \frac{(f_1-f_0)^2}{f_0} \dd\mu
  =
  1 + \chi^2(Q_0,Q_1),
\end{equation}
where $\chi^2$ stands for the $\chi^2$ distance between probability measures
(see, e.g., \cite[Sect.~31, Definition~2]{Borovkov:1998});
the value~\eqref{eq:after-changepoint} exceeds 1 unless $Q_0=Q_1$.
This suggests that the blue line is close to being horizontal before the changepoint
while it starts increasing after the changepoint.

While the calculations \eqref{eq:before-changepoint}--\eqref{eq:after-changepoint}
are useful and will be used in the proof of Proposition~\ref{prop:CUSUM-efficiency} below,
it would be wrong to interpret them directly as indicators of the tendency of the likelihood ratio martingale
to increase or decrease.
For that, we should use
\begin{equation*}
  \E_0 \ln L_n
  =
  \int f_0 \ln\frac{f_1}{f_0} \dd\mu
  =
  -\KL(Q_0,Q_1)
  <
  0
\end{equation*}
before the changepoint and
\begin{equation*}
  \E_1 \ln L_n
  =
  \int f_1 \ln\frac{f_1}{f_0} \dd\mu
  =
  \KL(Q_1,Q_0)
  >
  0
\end{equation*}
after the changepoint,
where $\KL$ stands for the Kullback--Leibler divergence.
Therefore, the slope of the blue line is expected
to be negative before the changepoint and positive after it.
This agrees with what we see in Figure~\ref{fig:Bernoulli}.

Let us now check informally that the conformal e-values $E_n$
and the likelihood ratios $L_n$ can be expected to be close to each other
soon after the changepoint.
Indeed, soon after the changepoint we will have
\begin{align}
  E_{N_0+n}
  &=
  \frac{L_{N_0+n}}{\frac{1}{N_0+n}(L_1+\dots+L_{N_0}+L_{N_0+1}+\dots+L_{N_0+n})}
  \label{eq:def}\\
  &\approx
  \frac{L_{N_0+n}}{\frac{N_0}{N_0+n}+\frac{1}{N_0+n}(L_{N_0+1}+\dots+L_{N_0+n})}
  \approx
  L_{N_0+n}.
  \label{eq:approx}
\end{align}
The equality in \eqref{eq:def} is just an application of the definition \eqref{eq:E-from-L}.
The first approximate equality (which holds with high probability) in \eqref{eq:approx}
follows from the law of large numbers assuming $N_0$ is large;
we obtained it by replacing $L_1+\dots+L_{N_0}$ by $N_0$
(see \eqref{eq:before-changepoint}).
In the denominator in \eqref{eq:approx} we have a weighted average of 1
and the postchange likelihood ratios.
The cumulative weight $\frac{n}{N_0+n}$ of the postchange likelihood ratios is small
if $n\ll N$,
and the second approximate equality in \eqref{eq:approx} assumes both $n\ll N$
and the postchange likelihood ratios being only moderately large.
The presence of the postchange likelihood ratios in the denominator in \eqref{eq:approx}
can be regarded as the usual conformal adjustment in this context.
The closeness of $E_n$ and $L_n$ is manifested in the closeness of the blue and orange lines
in Figure~\ref{fig:Bernoulli} soon after the changepoint.

The following proposition compares more formally
the behaviour of the likelihood ratio martingale
and the conformal e-pseudomartingale after the changepoint $N_0$.
For simplicity we will assume that the likelihood ratios $L_n$ are bounded;
e.g., $L_n\in[0.8,1.2]$ in the situation of Figure~\ref{fig:Bernoulli}.

\begin{proposition}\label{prop:CUSUM-efficiency}
  Suppose that $Z_n$, $n=1,\dots,N_0$, are generated from $Q_0\in\mathfrak{P}(\mathbf{Z})$,
  that $Z_n$, $n=N_0+1,\dots,N_0+N_1$, are generated from $Q_1\in\mathfrak{P}(\mathbf{Z})$,
  all independently,
  and that the likelihood ratios $L_n$ are bounded:
  $0\le a\le L_n\le b<\infty$.
  Suppose $\E_0 L_n=1$ and $c:=\E_1 L_n>1$.
  For any $\epsilon\in(0,1)$, we have
  \begin{multline}\label{eq:CUSUM-efficiency}
    \forall n\in\{1,\dots,N_1\}:
    \ln\frac{L_{N_0+1}\dots L_{N_0+n}}{E_{N_0+1}\dots E_{N_0+n}}
    <
    \frac{n(n+1)}{N_0}
    \frac{c}{2}\\
    +
    \frac{n}{N_0^{1/2}}
    \left|b-a\right|
    \sqrt{\frac12\ln\frac{2}{\epsilon}}
    +
    \frac{(N_1+1)^{3/2}}{N_0}
    \left|b-a\right|
    \sqrt{\frac16\ln\frac{2}{\epsilon}}
  \end{multline}
  with probability at least $1-\epsilon$.
\end{proposition}

The proof of Proposition~\ref{prop:CUSUM-efficiency}
is given in Appendix~\ref{app:proof-efficiency}.

The inequality in \eqref{eq:CUSUM-efficiency} shows that, for a fixed $\epsilon$,
the conformal e-pseudomartingale grows after the changepoint almost as fast
as the likelihood ratio martingale if $N_1\ll\sqrt{N_0}$.
In the situation of Figure~\ref{fig:Bernoulli},
the proposition tells us that we can expect similar rates of growth
for the conformal e-pseudomartingale and the likelihood ratio martingale
for around $\sqrt{1000}\approx30$ steps,
whereas we observe similar rates of growth for about $300$ steps.

\begin{figure}[bt]
  \begin{center}
    \includegraphics[width=0.48\textwidth]{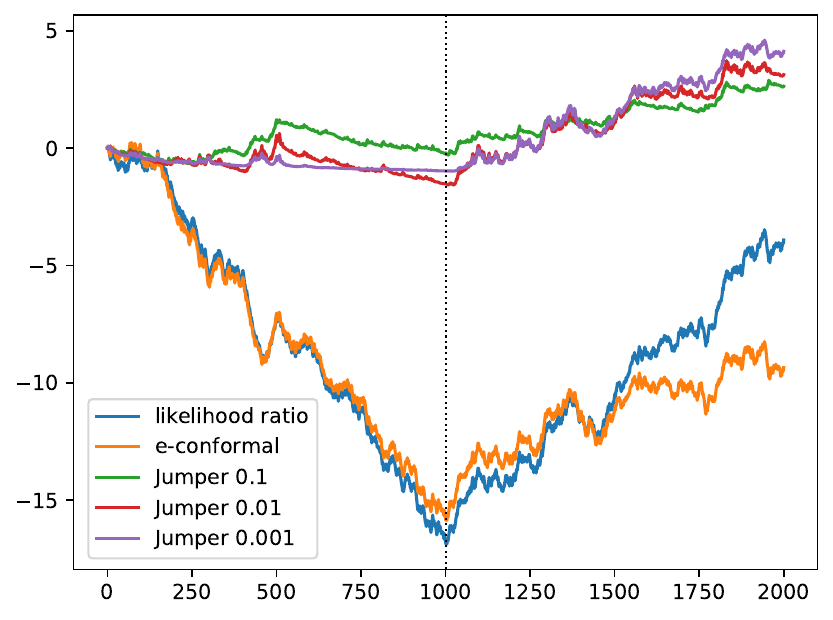}
    \includegraphics[width=0.48\textwidth]{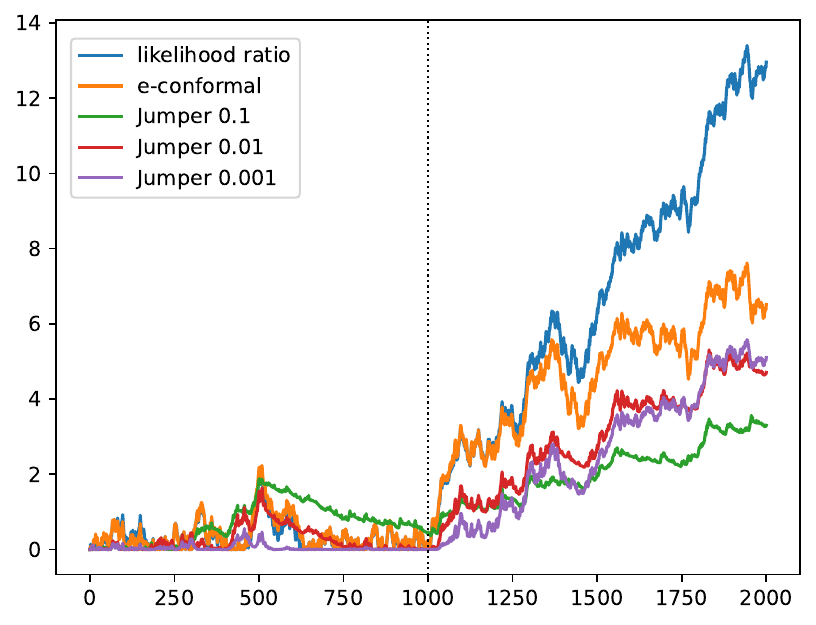}
  \end{center}
  \caption{Five processes as in Figure~\ref{fig:Bernoulli}
    but for the Cauchy distributions
    with the location and scale parameters $(0,1)$ before the changepoint
    and $(0,0.7)$ after the changepoint.
    Left panel: the raw processes.
    Right panel: the corresponding CUSUM statistics.}
  \label{fig:Cauchy}
\end{figure}

Figure~\ref{fig:Cauchy} is a counterpart of Figure~\ref{fig:Bernoulli}
for continuous, namely Cauchy, distributions.
Plots for Gaussian distributions, which are more standard,
look similar, but we have chosen a more awkward case of a distribution without a mean,
to make it less similar to the Bernoulli case.
The conformal test martingales, despite the generic nature of their betting functions,
now work much better than in Figure~\ref{fig:Bernoulli};
to see that they are still much worse than the conformal e-pseudomartingale,
the reader should pay attention to their behaviour soon after the changepoint at 1000,
where their growth is relatively sluggish.

\section{Conclusion}
\label{sec:conclusion}

The only known approach to detecting exchangeability violations online before this paper
was based on conformal prediction;
see, e.g., \cite[Part~III]{Vovk/etal:2022book}.
The approach of this paper is based instead on conformal e-prediction.
The two approaches are very different,
and neither dominates the other in all interesting applications.
These are some differences:
\begin{itemize}
\item
  Design of conformal test martingales involves two distinct steps:
  using a conformity measure to obtain p-values
  and then betting against those p-values.
  Conformal e-pseudomartingales do not involve such a rigid division
  and thus appear to be more flexible.
\item
  On the other hand, when betting on the $n$th step against the $n$th p-value $p_n$, $n=1,2,\dots$,
  conformal test martingales may use the previous p-values $p_1,\dots,p_{n-1}$.
  Such dependence on the past is not allowed for conformal e-pseudo\-mar\-tin\-gales.
\item
  Conformal test martingales are randomized
  (without randomization we only obtain conformal test supermartingales)
  whereas conformal e-pseudo\-mar\-tin\-gales do not require randomization
  (it is optional and not used in this paper).
\end{itemize}

We have discussed advantages and disadvantages of conformal testing and conformal e-testing.
In summary:
\begin{enumerate}
\item\label{it:1}
  In the online protocol,
  conformal testing relies on bona fide test martingales,
  which can then be used in one-off and multistage conformal testing.
  On the other hand, conformal e-pseudomartingales are not, in general,
  test martingales,
  and for them Ville's inequality can be violated badly
  (Sect.~\ref{sec:online}).
\item\label{it:2}
  A weakened variants of an efficiency result for conformal testing
  adapted to conformal e-testing
  is discussed in Sect.~\ref{sec:batch}.
\item\label{it:3}
  In Sect.~\ref{sec:CUSUM-validity} we construct a conformal CUSUM e-procedure
  that is a natural modification of the standard CUSUM procedure,
  while construction of an efficient conformal CUSUM procedure
  is more difficult and dependent on the postulated $Q_0$ and $Q_1$
  (whereas \eqref{eq:E-from-L} is applicable universally).
\item\label{it:4}
  A final advantage of conformal e-testing is that it requires no randomization.
\end{enumerate}
The first strength, \ref{it:1}, is a clear advantage
of conformal testing over conformal e-testing.
For strength \ref{it:2}, the picture is more ambiguous,
as the strength of conformal testing still partly survives for conformal e-testing.
And we also have two advantages, \ref{it:3} and \ref{it:4},
of conformal e-testing.

These are some interesting directions of further research:
\begin{itemize}
\item
  Developing betting functions that are better adapted
  to the assumed probability distributions $Q_0$ and $Q_1$
  than the generic ones used in the Simple Jumper method
  (Sect.~\ref{sec:CUSUM-efficiency}).
  (Perhaps in the spirit of the ``Bayes--Kelly'' approach of \cite[Part~III]{Vovk/etal:2022book}.)
\item
  Is it possible to strengthen Proposition~\ref{prop:CUSUM-efficiency}
  to obtain performance guarantees for the conformal CUSUM e-procedure
  that are comparable with what we observe in the experimental results?
\item
  One disadvantage of using the CUSUM procedure
  with a threshold $c$ and the target asymptotic frequency of false alarms $1/c$
  is that this method, while valid by Proposition~\ref{prop:CUSUM},
  might be a conservative way of achieving this target:
  namely, Proposition~\ref{prop:CUSUM} appears to be conservative,
  since even the Shiryaev--Roberts procedure achieves the target.
  Is it possible to employ, e.g.,
  adaptive conformal inference \cite{Gibbs/Candes:2021}
  or defensive forecasting \cite[Chap.~12]{Shafer/Vovk:2019}
  to adapt online the threshold $c$ to a target frequency of false alarms $\alpha$?
\end{itemize}

\subsection*{Acknowledgments}

This research has been partially supported by Astra Zeneca, Stena Line, and Mitie.

\appendix
\section{Proofs for Sect.~\ref{sec:CUSUM-validity}}
\label{app:proofs-validity}

The idea of the proof of Proposition~\ref{prop:rSR} consists,
as usual, in reversing the direction of time.
Let $N$ be a sufficiently large natural number
and $(\FFF_n)$ be the exchangeable filtration.

Remember that $E_n$ is the $n$th conformal e-variable \eqref{eq:E};
it has $e_n$ as it value.
Then $(E_n,\FFF_n)$, $n=N,\dots,1$, is an exact e-flow
(in the terminology of \cite{Vovk:sister})
over the finite time interval $N,\dots,1$,
in the sense
\[
  \E(E_n\mid\FFF_{n+1}) = 1,
  \quad
  n=N,\dots,1.
\]
The corresponding martingale is $(T_n,\FFF_n)$, $n=N+1,\dots,1$,
where
\[
  T_n
  :=
  E_n\dots E_N,
  \quad
  n=N+1,N,\dots,1,
\]
with $T_{N+1}$ understood to be $1$.
For simplicity,
let us assume that all $E_n$ are positive, so that $T_n$ is a positive martingale.

Let us apply the Shiryaev--Roberts procedure to the e-flow $E_{N},\dots,E_1,\dots$
continued by setting $E_0:=E_{-1}:=\dots:=1$.
It gives us the decreasing sequence of stopping times $\sigma_0:=N+1$ and
\begin{equation}\label{eq:SR}
  \sigma_k
  :=
  \max
  \left\{
    n<\sigma_{k-1}
    \mid
    \sum_{i=n}^{\sigma_{k-1}-1}
    E_n \dots E_i
    \ge
    c
  \right\},
  \quad
  k=1,2,\dots,
\end{equation}
where $n$ ranges over the integers.
A useful property of the Shiryaev--Roberts procedure is
$\E(\sigma_{k-1}-\sigma_k\mid\FFF_{\sigma_k})\ge c$;
see, e.g., \cite[Proposition~8.13]{Vovk/etal:2022book}
(proved in \cite[Sect.~8.5.2]{Vovk/etal:2022book};
the current proof also uses some other ideas in that section).

To make the stopping times \eqref{eq:SR} more manageable,
it will be convenient to force an alarm every $L$ steps,
where $L$ is to be chosen later:
$\sigma'_0:=N+1$ and
\begin{equation*}
  \sigma'_k
  :=
  (\sigma'_{k-1}-L)
  \vee
  \max
  \left\{
    n<\tau_{k-1}
    \mid
    \sum_{i=n}^{\sigma'_{k-1}-1}
    E_n \dots E_i
    \ge
    c
  \right\},
  \quad
  k=1,2,\dots;
\end{equation*}
let us call this the \emph{truncated Shiryaev--Roberts procedure}.
As in the proof of \cite[Corollary~8.4]{Vovk/etal:2022book},
induction in $k$ shows that $\sigma'_k\ge\sigma_k$ for all~$k$.

There is a useful connection between $\tau_k$ and $\sigma'_k$:
namely, each set $\{\tau_{k-1}+1,\dots,\tau_{k}\}$ with $\tau_k\le N$
contains at least one stopping time $\sigma'_l$.
(This is even true for $\sigma_l$, although we do not need it.)
This can be deduced from
\begin{equation}\label{eq:premise}
  \sum_{j=i}^{\tau_{k}}
  E_i\dots E_j
  \ge
  c
  \text{ for some $i\ge\tau_{k-1}+1$}.
\end{equation}
Indeed, let $l$ be the largest number satisfying $\sigma'_l>\tau_{k}$;
arguing indirectly, let us suppose that $\sigma'_{l+1}<\tau_{k-1}+1$.
The inequality \eqref{eq:premise} implies
\begin{equation*}
  \sum_{j=i}^{\sigma'_l-1}
  E_i\dots E_j
  \ge
  \sum_{j=i}^{\tau_{k}}
  E_i\dots E_j
  \ge
  c
  \text{ for some $i>\sigma'_{l+1}$},
\end{equation*}
which contradicts the definition of the stopping times $\sigma'$.

As in \cite[Proof of Proposition~8.15 in Sect.~8.5.2]{Vovk/etal:2022book},
let us say that $k$ is \emph{slow} if
\[
  \P
  \left(
    \sigma'_{k-1}-\sigma'_k = L
    \mid
    \FFF_{\sigma'_{k-1}}
  \right)
  \ge
  c/L
\]
and \emph{fast} otherwise.
We show there that
\begin{equation}\label{eq:EE}
  \E
  \left(
    \sigma'_{k-1}-\sigma'_k
    \mid
    \FFF_{\sigma'_{k-1}}
  \right)
  \ge
  c-c^2/L
\end{equation}
if $k$ is fast.
On the other hand, if $k$ is slow,
\[
  \E
  \left(
    \sigma'_{k-1}-\sigma'_k
    \mid
    \FFF_{\sigma'_{k-1}}
  \right)
  \ge
  L
  \P
  \left(
    \sigma'_{k-1}-\sigma'_k = L
    \mid
    \FFF_{\sigma'_{k-1}}
  \right)
  \ge
  c.
\]
In both cases, we have~\eqref{eq:EE}.

Let $A'_N$ be the largest $k$ such that $\sigma'_k>0$;
we interpret $A'_N$ as the number of alarms raised
by the truncated Shiryaev--Roberts procedure.
Remember that $A_N$, defined by \eqref{eq:A},
is the number of alarms raised by the conformal reverse Shiryaev--Roberts e\-/procedure.
By Hoeffding's inequality \cite[Sect.~A.6.3]{Vovk/etal:2022book},
for arbitrarily small $\epsilon\in(0,c)$ and sufficiently large $N$,
\begin{multline*}
  \P
  \left(
    \frac{A_N}{N}
    \ge
    \frac{1}{c-\epsilon}
  \right)
  \le
  \P
  \left(
    \frac{A'_N}{N}
    \ge
    \frac{1}{c-\epsilon}
  \right)
  =
  \P
  \left(
    A'_N\ge\left\lceil\frac{N}{c-\epsilon}\right\rceil
  \right)\\
  \le
  \P
  \left(
    \sum_{k=1}^{\lceil\frac{N}{c-\epsilon}\rceil}
    (\sigma'_{k-1}-\sigma'_k)
    \le
    N
  \right)
  \le
  \exp
  \left(
    -2
    \frac
      {\left(\lceil\frac{N}{c-\epsilon}\rceil\left(c-\frac{c^2}{L}\right)-N\right)^2}
      {L^2 \lceil\frac{N}{c-\epsilon}\rceil}
  \right)\\
  \le
  \exp(-\epsilon'N),
\end{multline*}
where $\epsilon'$ is a positive constant
(which requires a sufficiently large $L$).
Since the series $\sum_N\exp(-\epsilon'N)$ converges,
\[
  \frac{A_N}{N}
  \ge
  \frac{1}{c-\epsilon}
\]
happens only finitely often,
which completes the proof of Proposition~\ref{prop:rSR}.

To deduce Proposition~\ref{prop:CUSUM} notice that,
by induction in $k$,
$\tau_k^{\text{SR}}\le\tau_k^{\text{CUSUM}}$,
where $\tau_k^{\text{CUSUM}}$ are the $\tau_k$ defined by \eqref{eq:e-CUSUM},
and $\tau_k^{\text{SR}}$ are the $\tau_k$ defined by \eqref{eq:rSR}.

\section{Proof for Sect.~\ref{sec:CUSUM-efficiency}}
\label{app:proof-efficiency}

The proof will be based on the following ``maximal'' version of Hoeffding's inequality.

\begin{proposition}[maximal Hoeffding inequality]\label{prop:Hoeffding}
  Let $\FFF_0,\dots,\FFF_N$ be a filtration.
  For any deterministic sequence $c_1,\dots,c_N$ of positive numbers,
  any predictable sequence $a_1,\dots,a_N$ w.r.\ to $(\FFF_i)$,
  any supermartingale difference $\xi_1,\dots,\xi_N$ w.r.\ to $(\FFF_i)$
  such that $\xi_i\in[a_i,a_i+c_i]$, $i=1,\dots,N$,
  and any $\beta>0$,
  \begin{equation}\label{eq:H}
    \P
    \left\{
      \max_{n=0,\dots,N}
      \sum_{i=1}^n
      \xi_i
      \ge
      \beta
    \right\}
    \le
    \exp
    \left(
      -\frac{2\beta^2}{\sum_{i=1}^N c_i^2}
    \right)
  \end{equation}
\end{proposition}

For a proof of Proposition~\ref{prop:Hoeffding},
see, e.g., \cite[Theorem~3.2.1]{Roch:2024}
(or adapt a proof of Hoeffding's standard inequality,
such as that given in \cite[Sect.~A.6.3]{Vovk/etal:2022book}).

The interpretation of the three addends
on the right-hand side of the inequality in \eqref{eq:CUSUM-efficiency} is:
the first addend reflects the effect of the expectation $c$ of the postchange likelihood ratios,
the second addend reflects the volatility of the prechange likelihood ratios,
and the third addend reflects the volatility of the postchange likelihood ratios.
We start from the first addend and split the permitted probability $\epsilon$
of violating \eqref{eq:CUSUM-efficiency} into two equal parts,
one controlling the prechange behaviour of the likelihood ratios
and the other controlling their postchange behaviour.

By Hoeffding's inequality \eqref{eq:H} we will have
\begin{equation}\label{eq:prechange-bound}
  \sum_{n=1}^{N_0}
  L_n
  <
  N+\beta
  \text{ with probability}
  {}\ge 1-\frac{\epsilon}{2}
\end{equation}
when
\[
  \exp
  \left(
    -\frac{2\beta^2}{N_0(b-a)^2}
  \right)
  =
  \frac{\epsilon}{2},
\]
and so, solving this equation, we set 
\begin{equation}\label{eq:beta}
  \beta
  :=
  \sqrt{\frac12N_0(b-a)^2\ln\frac{2}{\epsilon}}.
\end{equation}

The log ratio in the left-hand side of the inequality in \eqref{eq:CUSUM-efficiency}
can be bounded from above as follows:
\begin{align}
  &\ln\frac{L_{N_0+1}\dots L_{N_0+n}}{E_{N_0+1}\dots E_{N_0+n}}\label{eq:log-ratio}\\
  &=
  \ln
  \left(
    \frac{L_1+\dots+L_{N_0}+L_{N_0+1}}{N_0+1}
    \dots
    \frac{L_1+\dots+L_{N_0}+L_{N_0+1}+\dots+L_{N_0+n}}{N_0+n}
  \right)\notag\\
  &<
  \ln
  \left(
    \frac{N_0+\beta+L_{N_0+1}}{N_0}
    \dots
    \frac{N_0+\beta+L_{N_0+1}+\dots+L_{N_0+n}}{N_0}
  \right)\notag\\
  &\le
  \frac{n\beta + n L_{N_0+1} +\dots+ L_{N_0+n}}{N_0},
  \label{eq:upper-bound}
\end{align}
where the first inequality, which holds with probability at least $1-\epsilon/2$,
follows from the inequality in \eqref{eq:prechange-bound}
(decreasing the denominators does not affect the validity of the first inequality),
and the second inequality follows from the standard inequality $\ln(1+x)\le x$.
The expectation of the upper bound \eqref{eq:upper-bound} on the log-ratio \eqref{eq:log-ratio} is
\begin{equation}\label{eq:expectation}
  \frac{n\beta + c(1+\dots+n)}{N_0}
  =
  \frac{n\beta + c n(n+1)/2}{N_0}.
\end{equation}
By Hoeffding's maximal inequality \eqref{eq:H} we will have
\begin{equation}\label{eq:postchange-bound}
  \forall n\in\{1,\dots,N_1\}:
  \frac{n(L_{N_0+1}-c)+\dots+(L_{N_0+n}-c)}{N_0}
  <
  B
\end{equation}
with probability at least $1-\epsilon/2$
when
\[
  \exp
  \left(
    -\frac{2B^2N_0^2}{(b-a)^2+\dots+N_1^2(b-a)^2}
  \right)
  =
  \frac{\epsilon}{2}.
\]
Solving the last equation and using the standard identity
$1+\dots+n^2=n(n+1)(2n+1)/6$,
we set
\begin{multline}\label{eq:B}
  B
  :=
  \frac{(N_1+1)^{3/2}}{N_0}
  \sqrt{\frac{1}{6}(b-a)^2\ln\frac{2}{\epsilon}}\\
  >
  \frac{1}{N_0}
  \sqrt{\frac{1}{12}N_1(N_1+1)(2N_1+1)(b-a)^2\ln\frac{2}{\epsilon}}.
\end{multline}

Combining \eqref{eq:beta}, \eqref{eq:log-ratio}--\eqref{eq:upper-bound},
\eqref{eq:expectation}, \eqref{eq:postchange-bound}, and \eqref{eq:B},
we obtain that \eqref{eq:CUSUM-efficiency} holds
with probability at least $1-\epsilon$,
which completes the proof of Proposition~\ref{prop:CUSUM-efficiency}.

\begin{thebibliography}{10}
\bibitem{Borovkov:1998}
Alexander~A. Borovkov.
\newblock {\em Mathematical Statistics}.
\newblock Gordon and Breach, Amsterdam, 1998.

\bibitem{Gibbs/Candes:2021}
Isaac Gibbs and Emmanuel~J. Cand\`es.
\newblock Adaptive conformal inference under distribution shift.
\newblock {\em Advances in Neural Information Processing Systems},
  34:1660–1672, 2021.

\bibitem{Lorden:1971}
Gary Lorden.
\newblock Procedures for reacting to a change in distribution.
\newblock {\em Annals of Mathematical Statistics}, 42:1897--1908, 1971.

\bibitem{Moustakides:1986}
George~V. Moustakides.
\newblock Optimal stopping times for detecting a change in distribution.
\newblock {\em Annals of Statistics}, 14:1379--1388, 1986.

\bibitem{Page:1954}
Ewan~S. Page.
\newblock Continuous inspection schemes.
\newblock {\em Biometrika}, 41:100--115, 1954.

\bibitem{Poor/Hadjiliadis:2009}
H.~Vincent Poor and Olympia Hadjiliadis.
\newblock {\em Quickest Detection}.
\newblock Cambridge University Press, Cambridge, 2009.

\bibitem{Ramdas/etal:2023}
Aaditya Ramdas, Peter Gr\"unwald, Vladimir Vovk, and Glenn Shafer.
\newblock Game-theoretic statistics and safe anytime-valid inference.
\newblock {\em Statistical Science}, 38:576--601, 2023.

\bibitem{Ritov:1990}
Ya'acov Ritov.
\newblock Decision theoretic optimality of the {CUSUM} procedure.
\newblock {\em Annals of Statistics}, 18:1464--1469, 1990.

\bibitem{Roch:2024}
S\'ebastien Roch.
\newblock {\em Modern Discrete Probability: An Essential Toolkit}.
\newblock Cambridge University Press, Cambridge, 2024.

\bibitem{Shafer/Vovk:2019}
Glenn Shafer and Vladimir Vovk.
\newblock {\em Game-The\-o\-ret\-ic Foundations for Probability and Finance}.
\newblock Wiley, Hoboken, NJ, 2019.

\bibitem{Shiryaev:2016}
Albert~N. Shiryaev.
\newblock {\em Probability-1}.
\newblock Springer, New York, third edition, 2016.

\bibitem{Ville:1939}
Jean Ville.
\newblock {\em Etude critique de la notion de collectif}.
\newblock Gauthier-Villars, Paris, 1939.

\bibitem{Vovk:2002FOCS}
Vladimir Vovk.
\newblock On-line {Confidence} {Machines} are well-calibrated.
\newblock In {\em Proceedings of the Forty Third Annual Symposium on
  Foundations of Computer Science}, pages 187--196, Los Alamitos, CA, 2002.
  IEEE Computer Society.

\bibitem{Vovk:2021-local}
Vladimir Vovk.
\newblock Testing randomness online.
\newblock {\em Statistical Science}, 36:595--611, 2021.
\newblock For earlier versions, see
  \href{https://arXiv.org/abs/1906.09256}{arXiv:\allowbreak 1906.09256}
  [math.PR] (from 2019).

\bibitem{Vovk:sister}
Vladimir Vovk.
\newblock Conformal e-prediction.
\newblock {\em Pattern Recognition}, 0:0--0, 2024.
\newblock Submitted for publication in the Special Issue on Conformal
  Prediction and Distribution-Free Uncertainty Quantification.

\bibitem{Vovk/etal:2022book}
Vladimir Vovk, Alex Gammerman, and Glenn Shafer.
\newblock {\em Algorithmic Learning in a Random World}.
\newblock Springer, Cham, second edition, 2022.

\bibitem{Vovk/etal:2003ICML}
Vladimir Vovk, Ilia Nouretdinov, and Alex Gammerman.
\newblock Testing exchangeability on-line.
\newblock In Tom Fawcett and Nina Mishra, editors, {\em Proceedings of the
  Twentieth International Conference on Machine Learning}, pages 768--775,
  Menlo Park, CA, 2003. AAAI Press.
\end{thebibliography}
\end{document}